\documentclass[11pt]{amsart}
\usepackage{amsfonts,amssymb,amsmath,amscd}

\usepackage[all]{xy}
\usepackage{url}
\usepackage{stmaryrd}
\usepackage{color}

\numberwithin{thm}{section}

\newcommand{\ga}[2]{\begin{gather}\label{#1}#2 \end{gather}}

\newcommand{\sO}{{\mathcal O}}

\newcommand{\C}{{\mathbb C}}

\newcommand{\F}{{\mathbb F}}

\newcommand{\N}{{\mathbb N}}
\renewcommand{\P}{{\mathbb P}}
\newcommand{\Q}{{\mathbb Q}}
\newcommand{\R}{{\mathbb R}}

\newcommand{\Z}{{\mathbb Z}}

\title{
Rational points over $C_1$ fields
}

\author{
  H\'el\`ene Esnault}
  \address{Freie Universit\"at Berlin, 
  Mathematik,
  Arnimallee 7,
  14195 Berlin, 
  Germany
  }
  \email{esnault@math.fu-berlin.de
   }

\begin{document}

\maketitle

\section{Introduction} \label{sec:intro}
Early in his work, Yuri Manin   established  a conjectural relation between    the geometry of certain smooth projective varieties and the  existence of rational points on them on some specific fields.  We shall focus on the first such example, published in 1966 in \cite{Man66}. He was 29. 

\medskip

In his PhD thesis ~\cite{Lan52},   Lang defines the notion of a quasi algebraically closed field, a concept he attributes to  E. Artin, see also the preface  \cite[p.x]{LT65} to Artin's Collected Papers by Lang and Tate:
 this is a field $k$ for which every hypersurface of degree $d\le n$ in $\P^n$ over $k$ admits a rational point.   He defines more generally a field to be $C_i$ for some natural number $i$ if every hypersurface of degree $d$ with $d^i\le n$ in $ \P^n$ over $k$ admits a rational point.   So quasi algebraically closed fields are precisely the $C_1$ fields.  He lists fields which are $C_1$: finite fields, according to the theorem of Chevalley (\cite[Th\'eor\`eme]{Che35}) and its refinement by  Warning   (see \cite[Satz~1a]{War35} in which he shows a  congruence modulo the characteristic  $p$ of $k$  for the number of rational points); function fields over an algebraically closed field, according to the theorem of 
 Tsen {\cite[Satz~5]{Tse34}  (who  had also developed a variant of the $C_i$-notion). Lang proves  in {\it loc. cit.} that complete discretely valued fields  with an algebraically  closed residue  field are $C_1$. In equal characteristic,   it includes the henselization  at a finite place of a function field over an algebraically closed field (as then rational points are dense in rational points over the completion), 
 in unequal characteristic it includes 
  finite extensions of the maximal unramified extension of $\Q_p$ (as again  rational points are dense in rational points over the completion).

\medskip

Else there are conjectures. For example  E. Artin conjectured that the maximal abelian extension $\Q^{\rm ab}$  of $\Q$ is $C_1$ (see the preface by Lang and Tate in {\it loc. cit.}).
More recently,  Fargues  conjectured in \cite[Conjecture~3.10]{Far20} that the field of functions $\bar \Q_p (FF)$  of the Fargues-Fontaine curve is $C_1$, in analogy  with Lang's conjecture on the function field  $\R(C_0)$  of the conic, $\ C_\circ={\rm Spec} \big(R[x,y]/(x^2+y^2+1)\big)$  over $\R$  without a point (\cite[p.~379]{Lan53}).

\medskip

Manin in {\it loc. cit.} Theorem~3.12 b)  proves  that del Pezzo surfaces $X$  over a finite field  $\F_q$ have a rational point. He uses the Weil conjectures for del Pezzo  surfaces: the second $\ell$-adic cohomology $H^2(X_{\bar \F_p}, \Q_\ell)$  is spanned by algebraic cycles, which are all defined over a finite extension of  $\F_q$, 
and the odd cohomology dies. So the Grothendieck-Lefschetz trace formula~\cite[(25)]{Gro64} expresses  the number of $\F_q$-points $|X(\F_q)|$ of $X$ in the form 
\ga{}{ \ (\star) \  \ \ |X(\F_q)|=  1+  q(\sum_{j=1}^{b_2}  \xi_j) + q^2  \notag} 
 where the $\xi_i$ are roots of unity.  This implies that $|X( \F_q)| \neq \emptyset$. In fact, more is true: 
 as $(\sum_{i=1}^{b_2} \xi_j)$ is an algebraic integer, it implies that  $|X(\F_q)|$ is congruent to $1$ modulo $q$. 
Manin refers to  Weil's article \cite[p.557]{Wei58} for the shape of the trace formula in $\ell$-adic cohomology.  Weil did not have $\ell$-adic cohomology at disposal, he was conjecturing its existence and purity, proved later by Deligne in \cite{Del74}, and in \cite{Del80}.  In {\it loc. cit.} Weil mentions  explicitly the formula $(\star)$ for rational surfaces  as a consequence of his conjecture. He argues that on those,  cohomological correspondences on $X\times_{\F_q} X$  are coming from algebraic cycles defined over some finite extension of $\F_q$. 
Even if he does not quote them, Manin likely had \cite{Gro64} at disposal, and    \cite{Tat65} 
in which Tate formulates the Tate conjecture and details the relation between the zeta function and the cohomology as a Galois module.  

\medskip

Manin  then formulates  in  {\it loc. cit.} Conjecture~4.1
   to the effect that every rational surface $X$  over a $C_1$ field  $k$ admits a rational point.   He proves it by  purely geometric methods  for  rational surfaces  with an extra geometric property, for example when the surface is fibered $f: X\to C$ in conics. Then   $C$ itself is a conic, that is a  curve of degree  $\le 2$  in $\P^2$, so has a rational point by the $C_1$ definition, and by the same argument its fibre as well.

 \medskip
 
 Manin's conjecture  was really  {\it programmatic}, which made arithmetic  and geometric properties  intertwined. It is interesting to mention that Deligne, to quote himself, coined the term {\it Arithmetic Geometry}, influenced by the title of the conference ``Arithmetic Algebraic Geometry''  held in Purdue in December 1963, in which Tate exposed his conjectures, but also after Manin's article {\it loc. cit.}.  Of course, the Weil conjecture, the Tate conjecture etc. postulate properties shared by all varieties over certain fields, while Manin initiates the study of the relation between some specific   fields and some specific varieties.
 
 \medskip

Campana in \cite{Cam91} in complex geometry and 
 Koll\'ar, Miyaoka and  Mori in \cite{KMM92} more generality develop the notion of 
 separably  rationally connected varieties, a vast generalization of rational varieties on one hand, and of hypersurfaces of degree $d\le n$ in $\P^n$ on the other.  In dimension $1$ and $2$, rationally connected varieties are precisely rational varieties.  
 In any dimension, a smooth projective connected variety $X$ over a field  $k$  is separably rationally connected if  
over the algebraic closure $\bar k$  it admits a  free rational curve, that is a morphism $ \P^1\to X_{\bar k}$ such that the pull-back of the tangent bundle is ample on the rational curve. If $k$ has characteristic $0$ this is equivalent to $X$ being  rationally  chain connected, that is to  the property that any two geometric points can be linked by a finite connected chain of rational curves.   If $k$ has characteristic $p>0$, if $X$ is separably rationally connected, it is  rationally  chain connected. 
  The positivity of the  pull-back of the tangent bundle forces those rational curves to be separable. 
 
 
 \medskip
 
 Koll\'ar publishes his book \cite{Kol96} in 1996 in which he studies general properties of rationally connected varieties. 
Manin's conjecture by this time has become the {\it  Lang-Manin conjecture}:
 separably rationally connected varieties should have a rational point over a $C_1$ field.  In the sequel we call it {\it the $C_1$-conjecture.} We mention also 
 \cite[6.1.1]{Kol96} in which Koll\'ar poses the conjecture as a problem over a function field and \cite[6.1.2]{Kol96} in which he asks for all fields over which separably  rationally connected varieties have a rational point.

\medskip

In this short note, we report on some progress on the  $C_1$-conjecture.   There are two main directions.

\medskip

 The first one relies on a classification of the varieties considered, see Section~\ref{sec:geom}. The proof by 
 Colliot-Th\'el\`ene~\cite[Proposition~2]{CT86} of Manin's initial conjecture, that is for rational surfaces, relies on the classification of rational surfaces over any field (by Iskovskikh ~\cite{Isk80} for the last stone). They  are fibered in conics,  the case 
already essentially treated by Manin, or del Pezzo surfaces,  see \cite[Th\'eor\`eme~1.11]{Wit10}.

\medskip

The second one consists in going through all the  known $C_1$ fields one by one and try to apply methods 
specific to those fields. For example if $k$ is a field of functions, do global geometry and degeneration on the special fibres, see Section~\ref{sec:GHS}.
If  $k$ is the henselization  at a finite place of a field of functions, reduce to the previous case, see 
Section~\ref{sec:hensel}.  If  $k$ is finite, make sure that even if $\ell$-adic cohomology is not spanned by algebraic cycles, the ``level'' of the cohomology is at least one so as to be able to apply $q$-divisibility in the Grothendieck-Lefschetz trace formula, see Section~\ref{sec:Esn}. If $k$  is a complete discretely valued field with an algebraically closed residue field in unequal characteristic, we know very little, see Section~\ref{sec:uneq} for a few modest  facts
 on the index and on a ``dynamical'' result.

\medskip 

On the other hand, we do not know a complete list of $C_1$ fields.  This shows that we are very far  from understanding the conjecture in all generality. At the same time,  there  are classical  varieties which are separably rationally connected and not defined by equations, but e.g. by a moduli functor, like in characteristic $0$  the moduli of stable bundles of coprime rank and degree on a curve (see \cite{Kau18} for first steps in this direction). It would be appealing if  over   a field  from which one conjectures  it is $C_1$, one could find a rational point on a separably rationally connected variety which is not defined by equations.  Or perhaps no point:  in the latter case, we would have to revise the hierarchy of problems!

\medskip

This small text is by no means exhaustive and selects a few spots which are closer to the author's interests and understanding.  For background and general studies, we  refer to \cite{Deb03}, \cite{Wit10} and \cite{CT11}.

\medskip

{\it Acknowledgements:} We thank Pierre Deligne for  a discussion on \cite{Wei58} and  on the origins of the 
terminology ``Arithmetic Geometry'' which is reflected in the Introduction.  We thank Lars Hesselholt for discussions on  $\bar \Q_p(FF)$, 
   Olivier Benoist and  Kay R\"ulling    for discussions on $\R(C_\circ)$,  Marco D'Addezio for discussions on   
   $\Q^{\rm ab}$, which are  reflected in Section~\ref{sec:brauer},  Olivier Wittenberg for general discussions on the topic when it came,  J\'anos Koll\'ar for precise references in his book, 
   Akhil  Mathew for his interest. 

\section{General Facts} \label{sec:coll}

\subsection{The Brauer group} \label{sec:brauer}
The Brauer group of a $C_1$ field is zero, see \cite[Corollary~of~Theorem~14]{Lan52}, see \cite[Proposition1.15]{Wit10} for a self-contained proof,  where the $C_1$ property is applied to the reduced norm of the  skew field to its center.

\medskip

In particular, to check whether a field  $k$ is $C_1$, one  might first compute  that its Brauer group  ${\rm Br}(k)$ is equal to zero.

\medskip

Witt proves in \cite[Satz~4]{Witt34}  that the norm map  is surjective on finite extensions of  $\R(C_\circ)$, which shows that ${\rm Br}( \R(C_\circ))=0$ (see \cite[Proposition~11, 3)]{Ser79}). 

\medskip

 For $\Q^{\rm ab}$, one proof is suggested in \cite[p.~162, d)]{Ser79}:  $\Q^{\rm ab}$ contains the maximal cyclotomic extension of $\Q$ (in fact it is equal to  it by the Kronecker-Weber theorem). The Brauer group of 
 a $p$-adic field   is killed by the maximal cyclotomic extension which is defined over $\Q$. The local-global principle  applied to the tower of cyclotomic extensions of $\Q$  
  allows one to conclude. 

\medskip
 Fargues proves $ {\rm Br}( \Q_p(FF))=0$ in  \cite[Th\'eor\`eme~2.2, Corollaire~2.5]{Far20}.

\subsection{Severi-Brauer varieties}
As the Brauer group of a  $C_1$ field $k$  is trivial,  Severi-Brauer varieties over $k$,
that is smooth projective varieties over $k$ which are isomorphic to $\P^d$ over an algebraic closure $\bar k$, 
 have a rational point.


\section{Geometry} \label{sec:geom}
In this section, we follow the Bourbaki talk by Debarre~\cite{Deb03}. 

\subsection{The theorem of Graber-Harris-Starr } \label{sec:GHS}
In \cite[Theorem~1.2]{GHS03} the authors prove the $C_1$ conjecture for  $k=\C(B)$  the field of functions of a curve $B$  over the field of complex numbers  $\C$.  This is equivalent to saying that if $f: X\to B$ is a projective morphism over a curve $B$ over $\C$,  with rationally connected generic fiber $ X\times_B \C(B)$, then $f$ has a section. 

\medskip

By replacing $f$ by its Weil restriction with respect to a Noether normalization $B\to \P^1$, we may assume that $B=\P^1$. 

\medskip

A closed point of $ X\times_{\P^1} \C(\P^1)$ corresponds to a diagram
\ga{}{ \xymatrix{\ar[d]_u  C \ar[r]^g & X \ar[dl]^f\\
\P^1} \notag}
where $C$ is a curve  and  $u: C\to \P^1$ is finite. If  $g$, viewed as a point in the 
moduli space of maps to $\P^1$ with certain properties at   the   ramification of $C\to g(C) $, was smooth,  then it would deform to a section of $f$ (\cite[Theorem~4.1]{Deb03}). This is the first place where  the authors use transcendental methods. The moduli $M(X, g_*[C])$  is the one of such maps for a fixed Betti class $g_*[C]\in H_2(X, \Z)$. They need its compactification for the proof. Then they have to put themselves in this situation. To this aim, they use topological arguments 
combined with the abundance of rational curves to prove that  once they have  $g$,  they can attach to $C$ 
 many rational curves  so the union deforms well. 

\subsection{The theorem of de Jong-Starr}  \label{sec:dJS}
In \cite[Theorem]{dJS03} the authors prove the $C_1$ conjecture for $k=F(B)$ where $F$ is any algebraically closed field. Strongly using the freeness of the rational curves, they deform directly algebraically the starting $g$ bypassing the use of $M(X, g_*[C])$ and its compactification, and the topological arguments.

\section{Motives and Cohomology} \label{sec:Esn}

\subsection{Grothendieck-Lefschetz trace formula}
As we saw in the introduction, if the whole $\ell$-adic cohomology of a smooth geometrically connected projective variety  $X$ is spanned by the cycle classes of algebraic cycles, then the Grothendieck-Lefschetz trace formula reads
\ga{}{ |X(\F_q)|= 1+  \sum_{i=1}^{2d} (-1)^j {\rm Tr}  \ {\rm Fr}| H^i(X_{\bar \F_p}, \Q_\ell)= 1+ \sum_{i=1}^d  q^i(\sum_{j_i=1}^{b_{2i}} \xi_{j_i})\notag}
where the $\xi_{j_i}$ are roots of unity,  the $b_{2i}$ are the Betti numbers of the even weighted  $\ell$-adic cohomology of $X_{\bar \F_p}$, ${\rm Fr}$ is the arithmetic Frobenius acting of $\ell$-adic cohomology.  (In addition the coefficient of $q^d$ is equal to $1$).  As $(\sum_{j=1}^{b_{2i}} \xi_{j})$ is an algebraic integer, 
from the formula  follows  the congruence
\ga{}{  |X(\F_q)| \equiv 1 \ {\rm mod} \ q . \notag}
 More generally, if the trace formula reads
\ga{}{  |X(\F_q)|= 1+ q \xi\notag}
where $\xi$ is an algebraic integer, then we can conclude that $ |X(\F_q)| \equiv 1 \ {\rm mod} \ q $. 

\medskip 

In \cite{Esn03},  we prove this property for rationally chain connected varieties, which  in particular proves the $C_1$ conjecture (without the  extra 
separability assumption) over finite fields. 

\medskip

We explain how  we reach this congruence in the sequel. In {\it loc. cit.} we used crystalline cohomology. The same argument is transposed here on the $\ell$-adic side. 

\subsection{Integrality: Deligne's theorem}
In \cite[Corollary~5.5.3]{Del73} Deligne proves that the eigenvalues of  ${\rm Fr}$ on $H^i(X_{\bar \F_p}, \Q_\ell)$ are algebraic integers. This is an important property in our context as it  shows:
\begin{quotation}
the  $C_1$ conjecture over finite fields follows from the statement that  the eigenvalues of ${\rm Fr}$ on $H^i(X_{\bar \F_p}, \Q_\ell)$ for $i>0$ are all divisible by $q$ as algebraic integers. 
\end{quotation}
\subsection{Motivic analogy} \label{sec:mot_ana}
In \cite[Introduction]{Esn06}, we raise the question of
\begin{quotation}  the analogue in $\ell$-adic cohomology    of the Hodge level  for  Betti cohomology. 
\end{quotation}
If $X$ is  a smooth  complex projective variety,  if $H^i(X, \Q)$ is supported in codimension $\ge 1$, then the Hodge level of $H^i(X, \Q)$  is $\ge 1$ as well.  The converse is an example of Grothendieck's generalized Hodge conjecture.  

\medskip

Using purity,  due to Gabber, we compute  in \cite[Lemma~2.1]{Esn03} that if $X$ is  a smooth   projective variety defined over $\F_q$, and  if $H^i(X_{\bar \F_p} , \Q_\ell)$ is supported in codimension $\ge 1$, then the eigenvalues of ${\rm Fr}$ on 
 $H^i(X_{\bar \F_p} , \Q_\ell)$ are divisible by $q$ as algebraic integers. 

\medskip

Thus for proving the $C_1$ conjecture over finite fields,  we are reduced to proving that $H^i(X_{\bar \F_p} , \Q_\ell)$ is supported in codimension $\ge 1$ for all $i\ge 1$. 

\medskip

In \cite[Theorem~1.3]{BER12} we relate  the Hodge level in characteristic $0$ to the analogue in crystalline cohomology  of  the $q$-divisibility of the eigenvalues of  ${\rm Fr}$  in $\ell$-adic cohomology: the former on a smooth projective variety in characteristic  $0$ implies  that the slopes in rigid cohomology $H^i$  of a  specialization of a regular model all are $\ge 1$ for $i\ge 1$ .  This makes in this context  the philosophical analogy a concrete theorem. 

\subsection{Trivial Chow group of $0$-cycles after base change to  any algebraic closed field}

By the very definition of rationally chain  connected varieties, any two geometric points are linked by a finite  connected chain of rational curves. In particular, the group  $CH_0(X_K)$ of algebraic $0$-cycles is equal to $\Z$  for any embedding   $k \hookrightarrow K$ of $k$ into an algebraically closed field $K$. Taking $K$ to contain the generic point $ {\rm Spec}(k(X))$,  the resulting decomposition of the diagonal initiated by  Bloch in his book~\cite[Appendix~to~Lecture~1]{Blo80}, a method which later became the basis of manifold developments, 
proves that $H^i(X_{\bar \F_p} , \Q_\ell)$ is supported in codimension $\ge 1$ for all $i\ge 1$.  This finishes the proof of $C_1$ conjecture over finite fields.

\section{Remarks} \label{sec:rmks}

\subsection{ $C_1$ conjecture and minimal model program in characteristic $0$}
A generalization of the study   case by case of rational surfaces by Colliot-Th\'el\`ene {\it loc. cit.} is performed in characteristic $0$  by Pieropan  in  \cite{Pie22}: over a given field of characteristic $0$, the minimal model program implies that 
all rationally connected varieties of dimension $\le d$ have a rational point  if and only if $\Q$-factorial Fano varieties of dimension $\le d$ and Picard rank $1$ do. It is roughly  a generalization of the reduction from all rational surfaces to the del Pezzo ones. 

\subsection{Henselization of function fields at finite places} \label{sec:hensel}
In \\  \cite[Th\'eor\`eme~7.5]{CT11},  Colliot-Th\'el\`ene proves that 
the $C_1$ conjecture holds true for   $k$ being  the field of fractions of an  henselian discrete  valuation ring 
with algebraically closed residue field $F$ in equal characteristic.  Indeed, if $X$ is defined over $k$, its rational points $X(k)$ are dense  in $X(k^{\widehat{}}) $ for the topology defined by the discrete valuation, where $k^{\widehat{} }=F((t))$   is the completion with respect to the valuation. If $X$ was defined over $k(t)$, we could then apply directly  the theorems   Sections~\ref{sec:GHS} and ~\ref{sec:dJS}. In general $X$ descends to $X_A$ where  $A=F[T_1,\ldots, T_n]/I \subset F[[t]]$ is an affine  and  smooth algebra over $F$, where the embedding is defined by a power series expansion in $t$ of the $T_i$.  For $n \in \N_{\ge 1} $ large enough, depending on the multiplicity of ${\rm Spec} F\in S={\rm Spec}(A)$,  the $F[[t]]/(t^n)$-point of $S$ integrates to a finite type normal curve $C_n \to S$  over $F$, generically embedded in $S$,  thus over which $X_{C_n}$ has a $C_n$-point.  In particular, $X(F[[t]]/(t^n))\neq \emptyset.$  As $X$ is proper, 
a compactness argument allows one to conclude that $X(F[[t]])\neq \emptyset$.

\subsection{Some remarks on complete discrete valuation rings $\sO$  with algebraically closed residue field $F$ in unequal characteristic } \label{sec:uneq}
\subsubsection{ A small dream}
We mention now where the  proof from Section~\ref{sec:hensel}  can not be generalized to the unequal characteristic situation. By the same argument as in  Section~\ref{sec:GHS}, we may assume that $\sO=W(F)$, where $F$ is the  algebracially closed residue field of characteristic $p>0$.  So we have $A\subset W$, which is now defined by the $p$-adic expansion of the $T_i$,  the resulting $W_n=W/p^n$ points of $S$,
$X_A$.  But we do not have the integration $C_n$ as in the equal characteristic case.  At the ICM 2014 in Seoul, we asked Scholze, which had just proved  Deligne's  weight-monodromy conjecture \cite[Theorem~1.14]{Sch12} for complete intersections  by applying his tilting method, whether he would think there is a way to tilt $X_W$    to  say $Y_{F[[t]]}$, keeping the rational connectivity property for $Y_{F((t))}$, and if so if one could   tilt back
 $Y_{F[[t]]}(F[[t]])$ constructed  in Section~\ref{sec:hensel} 
 to $X_W(W)$.  It is a difficult and badly posed question for many reasons, two immediate ones being that the tilt of  the field of fractions of $W$ contains also all the $p$-roots of $t$, and that  over this large field, $Y$ is  not rationally connected.  Nonetheless we can dream of a more evolved analogue of the classical proof in Section~\ref{sec:hensel} in this context. 

\subsubsection{Back to the Hodge level}  In \cite[Corollary~3]{ELW15} a weak form of the $C_1$ conjecture is proved over $K$ the field of fractions of $\sO$:   the index of a separably  rationally  connected variety $X$,  that is the g.c.d. of the degrees of  its closed points, is equal to $1$ if ${\rm dim} X +1<p$.   It is really a very weak form. Indeed, the index $1$ statement is true as soon as  the Hodge level of $X$ is $\ge 1$.
It is in the spirit of the theorem~\cite[Theorem~1.3]{BER12} mentioned in Section~\ref{sec:mot_ana}

\subsubsection{Dynamic}
As explained previously, fixing a  field $K$, complete  with respect to a discrete valuation with  algebraically closed residue field in unequal characteristic, the $C_1$ conjecture is not understood. In \cite[Theorem~1.3]{DK17}, Duesler and Knecht prove a ``dynamical'' version of the conjecture which is close to the theorem of  Ax-Knochen according to which, except for finitely many $p$,  any hypersurface   in $\P^n$ over $\Q_p$ of degree $d$  with $d^2\le n$ has a rational point. 
 They fix a Hilbert polynomial $P$ and show using similar methods  that, except for finitely many $p$,
 any  rationally connected variety with Hilbert polynomial $P$ and defined over the maximal unramified extension of $\Q_p$ has a rational point.


\begin{thebibliography}{DK09-2}

\bibitem[BER12]{BER12} Berthelot, P., Esnault, H., R\"ulling, K.: {\it Rational points over finite fields for regular models of algebraic varieties of Hodge type greater or equal to $1$},  
Annals of Mathematics  {\bf 176}  (2012), 1--96. 


\bibitem[Blo80]{Blo80} Bloch,S.: {\it Lectures on Algebraic cycles}, Duke University Mathematics Series {\bf IV} (1980).
\bibitem[Cam91]{Cam91} Campana, F.: {\it On twistor spaces of the class $C$}, J. Diff. Geom. {\bf 33} (1991), 541--549.

\bibitem[Che35]{Che35} Chevalley, C. {it D\'emonstration d'une hypoth\`ese de M. Artin}, Hamburg Abh. {\bf 11} (1935), 73--75.

\bibitem[CT86]{CT86} Colliot-Th\'el\`ene, J.-L.: {\it Arithm\'etique des vari\'et\'es rationnelles et probl\`emes birationnels},  Proceedings  ICM 1986, vol. 1, Berkeley, Amer. Math. Soc., (1987), 641--653.


\bibitem[CT11]{CT11} Colliot-Th\'el\`ene, J.-L.: {\it Vari\'et\'es  presque rationnelles, leurs points rationnels et leurs d\'eg\'en\'erescences},  In Arithmetic Geometry,  Lecture Notes in Math. {\bf 2009} (2011),  1--44. Springer Verlag.


\bibitem[Deb03]{Deb03} Debarre, O.: {\it Vari\'et\'es rationnellement connexes}, Expos\'e  {\bf 905}, S\'eminaire Bourbaki 2001/2002, Ast\'erisque {\bf 290} (2003), 243--266.

\bibitem[DK17]{DK17} Duesler, B., Knecht, A.: {\it Rationally connected varieties over the maximally unramified extension of $p$-adic fields}, Rocky Mountain J. Mat. {\bf 47} (2017), (8), 2605--2617.

\bibitem[dJS03]{dJS03} de Jong, J., Starr, J.: {\it Every rationally connected variety over the function field of a curve has a rational point}, Amer. J. Math. {\bf 125} (2003) (3), 567--580.

\bibitem[Del73]{Del73} Deligne, P.: {\it  Th\'eor\`eme d’int\'egralit\'e},  Appendix to 
``Le eniveau de la cohomologie des intersections compl\`etes''  by N. Katz, Expos'e XXI in SGA 7, Lecture Notes in Math. {\bf 340}, 363--400, Springer-Verlag, New York, 1973.

\bibitem[Del74]{Del74} Deligne, P.: {\it La conjecture de Weil I}, Publ. math. IHES {\bf 43} (1974), 273--308.


\bibitem[Del80]{Del80} Deligne, P.: {\it La conjecture de Weil II}, Publ. math. IHES {\bf 52} (1974), 137--252.


\bibitem[Esn03]{Esn03} Esnault, H.: {\it  Varieties over a finite field with trivial Chow group of $0$-cycles have a rational point}, 
Invent. math. {\bf 151} (2003), 187--191. 


\bibitem[Esn06]{Esn06} Esnault, H.: {\it Deligne's integrality theorem in unequal characteristic and rational points over finite fields}, 
Annals of Mathematics  {\bf 164} (2006), 715--730. 


\bibitem[ELW15]{ELW15} Esnault, H., Levine, M., Wittenberg, O.: {\it Index of varieties over Henselian fields and Euler characteristic of coherent sheaves},  
J. Algebraic Geometry {\bf 24}  (2015), 693--718. 

\bibitem[Far20]{Far20} Fargues, L. : {\it  $G$-torseurs en th\'eorie de Hodge $p$-adique},  Compositio Math.  {\bf 156} (2020), 2076--2110. 

\bibitem[GHS03]{GHS03} Graber, T., Harris, J., Starr, J.: {\it Families of ratonally connected varieties}, J. AMS {\bf 16} (2003), (1), 57--67. 

\bibitem[Gro64]{Gro64} Grothendieck, A.: {\it  Formule de Lefschetz et rationalit\'e  des fonctions $L$}, S\'em. Bourbaki {\bf 279}, 41--55, Soc. Math. France, Paris (1964/1965), 1995.

\bibitem[Isk80]{Isk80} Iskovskikh, V.: {\it Minimal models of rational surfaces over arbitrary fields}, Izv. Akad. Nauk SSSR Ser. Mat. {\bf 43} (1979), no. 1, 19--43 (in russian),   Math. USSR Izv. {\bf 14} (1980), no. 1.
(in english).

\bibitem[Kau18]{Kau18} Kaur, I.: {\it A pathological case of the $C_1$ conjecture in mixed characteristic}, In Mathematical Proceedings of the
Cambridge Philosophical Society,  {\bf 167} (2019), 61--64. 



\bibitem[KMM92]{KMM92} Koll\'ar, J., Miyaokoa, Y., Mori, S.: {\it Rational connectedness and boundedness of Fano manifolds}, J. Diff. Geom. {\bf 36} (1992) 765--779. 

\bibitem[Kol96]{Kol96} Koll\'ar, J.: {\it Rational Curves on Algebraic  Varieties}, Ergebnisse der Mathematik und ihrer Grenzgebiete {\bf 32}, Springer Verlag (1996).

\bibitem[Lan52]{Lan52} Lang, S.: {\it On quasi-algebraic closure}, Annals of Math. {\bf 55} (1952), 373--390. 

\bibitem[Lan53]{Lan53} Lang, S.: {\it The theory of real places}, Annals of Math. {\bf 57} (1957) (2), 378--391.

\bibitem[LT65]{LT65}  Lang, S., Tate, J.: {\it The collected papers of Emil Artin}, Addison-Wesley Publishing Co. (1965), xvi+560 pp.

\bibitem[Man66]{Man66} Manin, Yu.: {\it Rational surfaces over perfect fields},  Publ. math. 
IHES  {\bf 30} (1966), 55--113.


\bibitem[Pie22]{Pie22} Pieropan, M. {\it On rationally connected varieties over $C_1$ fields of characteristic $0$}, Algebra Number Theory {bf 16} (2022), (8), 1811--1844. 


\bibitem[Sch12]{Sch12} Scholze, P.: {\it  Perfectoid a spaces}, Publ. math. IHES {\bf 116} (2012), (1), 245--313.

\bibitem[Ser79]{Ser79} Serre, J.-P.:  {\it  Local Fields}, Graduate Texts in Mathematics {\bf 67}, Springer Verlag  (1979).

\bibitem[Tat65]{Tat65} Tate, J.: {\it Algebraic cycles and poles of zeta functions}, in {\it Arithmetical Algebraic Geometry},  Proc. Conf. Purdue Univ., (1963) pp. 93--110 Harper and  Row, New York (1965).

\bibitem[Tse34]{Tse34} Tsen, Ch. C.: {\it Algebren \"uber Funktionenk\"orpern}, Ph. D. dissertation, G\"ottingen University (1934).



\bibitem[War35]{War35} Warning, E.: {\it Bemerkung zur vorstehender Arbeit von Herrn Chevalley}, Hamburg Abh. {\bf 11} (1935), 76--83.


\bibitem[Wei58]{Wei58} Weil, A.: {\it Abstract versus classical algebraic geometry}, Matematika {\bf 4}  Vol.  2,  59--66.

\bibitem[Witt34]{Witt34} Witt, E.: {\it Zerlegung reeller algebraischer Funktionen in Quadrate. Schiefk\"orper \"uber reellem Funktionenk\"orper}, 
J. Reine Angew. Math. {\bf 171} (1934), 4--11.

\bibitem[Wit10]{Wit10} Wittenberg, O.: {\it La connexit\'e rationnelle en arithm\'etique}, in {\it Vari\'et\'es rationnellement connexes: aspects g\'eom\'etriques et arithm\'etiques}, Panorama et Synth\`eses {\bf 31}, Soci\'et\'e math\'ematique de France (2010). 















 



\end{thebibliography}
\end{document}